\documentclass{amsart}

\usepackage{amsmath,amsthm,amsfonts,amscd,amssymb}
\usepackage[pdfborder={0 0 0}]{hyperref}
\usepackage{enumerate}

\newcommand{\ZZ}{\mathbb{Z}}

\newtheorem{thm}{Theorem}[section]

\newtheorem{lem}[thm]{Lemma}
\newtheorem{prop}[thm]{Proposition}

\newtheorem{quest}{Question}

\theoremstyle{definition}

\theoremstyle{remark}
\newtheorem{rem}{Remark}[section]

\begin{document}

\title{On integral well-rounded lattices in the plane}
\author[L. Fukshansky]{Lenny Fukshansky}
\author[G. Henshaw]{Glenn Henshaw}
\author[P. Liao]{Philip Liao} 
\author[M. Prince]{Matthew Prince}
\author[X. Sun]{Xun Sun}
\author[S. Whitehead]{Samuel Whitehead}\thanks{The authors were supported by a grant from the Fletcher Jones Foundation. The first author was also partially supported by a grant from the Simons Foundation (\#208969 to Lenny Fukshansky) and by the NSA Young Investigator Grant \#1210223.}

\address{Department of Mathematics, 850 Columbia Avenue, Claremont McKenna College, Claremont, CA 91711}
\email{lenny@cmc.edu}
\address{Department of Mathematics and Computer Science, Wesleyan University, Middletown, CT 06459}
\email{ghenshaw@wesleyan.edu}
\address{Department of Mathematics, Claremont McKenna College, Claremont, CA 91711}
\email{PLiao14@students.claremontmckenna.edu}
\address{Department of Mathematics, Harvey Mudd College, Claremont, CA 91711}
\email{mthwate@gmail.com}
\address{School of Mathematical Sciences, Claremont Graduate University, Claremont, CA 91711}
\email{foxfur\_32@hotmail.com}
\address{Department of Mathematics, Pomona College, Claremont, CA 91711}
\email{scw22009@mymail.pomona.edu}


\begin{abstract} 
We investigate distribution of integral well-rounded lattices in the plane, parameterizing the set of their similarity classes by solutions of the family of Pell-type Diophantine equations of the form $x^2+Dy^2=z^2$ where $D>0$ is squarefree. We apply this parameterization to the study of the greatest minimal norm and the highest signal-to-noise ratio on the set of such lattices with fixed determinant, also estimating cardinality of these sets (up to rotation and reflection) for each determinant value. This investigation extends previous work of the first author in the specific cases of integer and hexagonal lattices and is motivated by the importance of integral well-rounded lattices for discrete optimization problems. We briefly discuss an application of our results to planar lattice transmitter networks.
\end{abstract}

\maketitle

\def\A{{\mathcal A}}
\def\AA{{\mathfrak A}}
\def\B{{\mathcal B}}
\def\C{{\mathcal C}}
\def\D{{\mathcal D}}
\def\EE{{\mathfrak E}}
\def\F{{\mathcal F}}
\def\x{{\mathcal H}}
\def\I{{\mathcal I}}
\def\II{{\mathfrak I}}
\def\J{{\mathcal J}}
\def\K{{\mathcal K}}
\def\kk{{\mathfrak K}}
\def\L{{\mathcal L}}
\def\LL{{\mathfrak L}}
\def\M{{\mathcal M}}
\def\mm{{\mathfrak m}}
\def\MM{{\mathfrak M}}
\def\N{{\mathcal N}}
\def\O{{\mathcal O}}
\def\OO{{\mathfrak O}}
\def\PP{{\mathfrak P}}
\def\R{{\mathcal R}}
\def\PNR{{\mathcal P_N(\real)}}
\def\PMNR{{\mathcal P^M_N(\real)}}
\def\PdNR{{\mathcal P^d_N(\real)}}
\def\s{{\mathcal S}}
\def\V{{\mathcal V}}
\def\X{{\mathcal X}}
\def\Y{{\mathcal Y}}
\def\Z{{\mathcal Z}}
\def\H{{\mathcal H}}
\def\cee{{\mathbb C}}
\def\Nn{{\mathbb N}}
\def\pee{{\mathbb P}}
\def\que{{\mathbb Q}}
\def\QQ{{\mathbb Q}}
\def\real{{\mathbb R}}
\def\RR{{\mathbb R}}
\def\zed{{\mathbb Z}}
\def\ZZ{{\mathbb Z}}
\def\aaa{{\mathbb A}}
\def\ff{{\mathbb F}}
\def\HDelta{{\it \Delta}}
\def\kk{{\mathfrak K}}
\def\qbar{{\overline{\mathbb Q}}}
\def\kbar{{\overline{K}}}
\def\ybar{{\overline{Y}}}
\def\kkbar{{\overline{\mathfrak K}}}
\def\ubar{{\overline{U}}}
\def\eps{{\varepsilon}}
\def\ahat{{\hat \alpha}}
\def\bhat{{\hat \beta}}
\def\gt{{\tilde \gamma}}
\def\h{{\tfrac12}}
\def\be{{\boldsymbol e}}
\def\bei{{\boldsymbol e_i}}
\def\bc{{\boldsymbol c}}
\def\bm{{\boldsymbol m}}
\def\bk{{\boldsymbol k}}
\def\bi{{\boldsymbol i}}
\def\bl{{\boldsymbol l}}
\def\bq{{\boldsymbol q}}
\def\bu{{\boldsymbol u}}
\def\bt{{\boldsymbol t}}
\def\bs{{\boldsymbol s}}
\def\bv{{\boldsymbol v}}
\def\bw{{\boldsymbol w}}
\def\bx{{\boldsymbol x}}
\def\bX{{\boldsymbol X}}
\def\bz{{\boldsymbol z}}
\def\bwy{{\boldsymbol y}}
\def\bY{{\boldsymbol Y}}
\def\bL{{\boldsymbol L}}
\def\ba{{\boldsymbol a}}
\def\bb{{\boldsymbol b}}
\def\bet{{\boldsymbol\eta}}
\def\bxi{{\boldsymbol\xi}}
\def\bo{{\boldsymbol 0}}
\def\bone{{\boldsymbol 1}}
\def\bol{{\boldsymbol 1}_L}
\def\ep{\varepsilon}
\def\p{\boldsymbol\varphi}
\def\q{\boldsymbol\psi}
\def\rank{\operatorname{rank}}
\def\aut{\operatorname{Aut}}
\def\lcm{\operatorname{lcm}}
\def\sgn{\operatorname{sgn}}
\def\spn{\operatorname{span}}
\def\md{\operatorname{mod}}
\def\Norm{\operatorname{Norm}}
\def\dim{\operatorname{dim}}
\def\det{\operatorname{det}}
\def\Vol{\operatorname{Vol}}
\def\rk{\operatorname{rk}}
\def\ord{\operatorname{ord}}
\def\ker{\operatorname{ker}}
\def\div{\operatorname{div}}
\def\Gal{\operatorname{Gal}}
\def\GL{\operatorname{GL}}
\def\SNR{\operatorname{SNR}}
\def\WR{\operatorname{WR}}
\def\IWR{\operatorname{IWR}}
\def\scg{\operatorname{\left< \Gamma \right>}}
\def\swrh{\operatorname{Sim_{WR}(\Lambda_h)}}
\def\ch{\operatorname{C_h}}
\def\cht{\operatorname{C_h(\theta)}}
\def\scgt{\operatorname{\left< \Gamma_{\theta} \right>}}
\def\scgmn{\operatorname{\left< \Gamma_{m,n} \right>}}
\def\gat{\operatorname{\Omega_{\theta}}}
\def\mn{\operatorname{mn}}
\def\disc{\operatorname{disc}}

\section{Introduction and statement of results}
\label{intro}

Let $N \geq 1$ be an integer, and let $\Lambda \subset \real^N$ be a lattice of full rank. Given a basis $\ba_1,\dots,\ba_N$ for $\Lambda$ we can write $A=(\ba_1 \dots \ba_N)$ for the corresponding basis matrix, and then $\Lambda = A \zed^N$. The corresponding norm form is defined as 
$$Q_A(\bx) = \bx^t A^t A \bx,$$
and we say that the lattice is {\it integral} if the coefficient matrix $A^t A$ of this quadratic form has integer entries; it is easy to see that this definition does not depend on the choice of a basis. The matrix $A^tA$ is called a {\it Gram matrix} of the lattice $\Lambda$. Integral lattices are central objects in arithmetic theory of quadratic forms and in lattice theory. We define $\det(\Lambda)$ to be $|\det(A)|$, again independent of the basis choice, and (squared) {\it minimum} or {\it minimal norm}
$$|\Lambda| = \min \{ \|\bx\|^2 : \bx \in \Lambda \setminus \{\bo\} \} = \min \{ Q_A(\bwy) : \bwy \in \zed^N \setminus \{\bo\} \},$$
where $\|\ \|$ stands for the usual Euclidean norm. Then each $\bx \in \Lambda$ such that $\|\bx\|^2 = |\Lambda|$ is called a {\it minimal vector}, and the set of minimal vectors of $\Lambda$ is denoted by $S(\Lambda)$. A lattice $\Lambda$ is called {\it well-rounded} (abbreviated WR) if the set $S(\Lambda)$ contains $N$ linearly independent vectors. These vectors do not necessarily form a basis for lattices in any dimension $N$, however they are known to form a basis for all $N \leq 4$ (see, for instance~\cite{pohst}); we will refer to such a basis as a {\it minimal basis} for $\Lambda$. WR lattices are important in discrete optimization, in particular in the investigation of sphere packing, sphere covering, and kissing number problems (see \cite{martinet}), as well as in coding theory (see \cite{esm}). Properties of WR lattices have also been investigated in \cite{mcmullen} in connection with Minkowski's conjecture and in \cite{lf:robins} in connection with the linear Diophantine problem of Frobenius. A particularly interesting and important class of WR lattices are the integral well-rounded lattices (abbreviated IWR). The main objective of the current paper is to study the properties of IWR lattices in the plane, extending some of the previous results of ~\cite{wr1}, \cite{wr2}, and \cite{wr3} with a view toward discrete optimization problems. Specifically, our investigation is motivated by the following three questions, which are the direct analogues of the questions asked about sublattices of the hexagonal lattice in~\cite{sloane}.

\begin{quest} \label{quest2} Which IWR lattice $\Lambda$ of a fixed determinant $\Delta$ maximizes the minimal norm? Since the density of circle packing associated to $\Lambda$ is equal to $\pi |\Lambda|/\Delta$, this choice of $\Lambda$ also maximizes the packing density.
\end{quest}

\begin{quest} \label{quest3} Which IWR lattice $\Lambda$ of a fixed determinant $\Delta$ maximizes the signal-to-noise ratio (defined below)? 
\end{quest}

\begin{quest} \label{quest1} How many IWR lattices of a fixed determinant $\Delta$ are there, up to rotation and reflection? This number is known to be finite.
\end{quest}

Given a lattice $\Lambda \in \real^N$, we can regard its nonzero points as transmitters which interfere with the transmitter at the origin, and then a standard measure of the {\it total interference} of $\Lambda$ is given by $E_{\Lambda}(2)$, where
\begin{equation}
\label{epstein_z}
E_{\Lambda}(s) = \sum_{\bx \in \Lambda \setminus \{ \bo \}} \frac{1}{\|\bx\|^{2s}}
\end{equation}
is the Epstein zeta-function of $\Lambda$, and the signal-to-noise ratio of $\Lambda$ is defined by
\begin{equation}
\label{SNR_def}
\SNR(\Lambda) = 10 \log_{10} \frac{1}{9E_{\Lambda}(2)},
\end{equation}
as in \cite{sloane}. To maximize $\SNR(\Lambda)$ on the set of all planar IWR lattices of a fixed determinant $\Delta$ is the same as to minimize $E_{\Lambda}(2)$ on this set. In fact, $E_{\Lambda}(s)$  for each real $s \geq 3$ is maximized by the same planar WR lattice of fixed determinant $\Delta$ that maximizes $|\Lambda|$, and vice versa (this follows from an old result of S. S. Ryskov \cite{ryskov}; see Lemma~\ref{min_epstein} below). Moreover, Lemma~\ref{SNR_min} and Remark~\ref{E_min_heuristics} below suggest that it may likely be so for $s=2$ as well. This would mean that Questions~\ref{quest2} and~\ref{quest3} are equivalent, which is not always so for non-WR lattices, as demonstrated in \cite{sloane}. 

Suppose that we have a network of transmitters positioned at the points of a planar lattice $\Lambda$. The plane is tiled with translates of the Voronoi cell of $\Lambda$, which are the cells serviced by the corresponding transmitters at their centers. The packing density of $\Lambda$ is precisely the proportion of the plane covered by the transmitter network. IWR lattices allow for transmitters of the same power and for integral distances between transmitters in the network, which simplifies positioning. Hence a lattice that answers Question~\ref{quest2}  maximizes coverage and the lattice that answers Question~\ref{quest3} maximizes signal-to-noise ratio for a 2-dimensional lattice transmitter network with a fixed cell area $\Delta$ and integral distances between transmitters. Our Lemma~\ref{SNR_min} and Remark~\ref{E_min_heuristics} below suggest that this may be done simultaneously. We present an algorithm for finding a lattice answering Question~\ref{quest2} for each possible value of $\Delta$ in Theorem~\ref{optimize}. It is also interesting to understand how many choices for positioning a network of equal-power transmitters with integral distances between them and fixed cell area are there -- this is an application of Question~\ref{quest1}; we estimate this number in Theorem~\ref{count}. We refer the reader to~\cite{heuvel} for further information about transmitter networks on planar lattices.
\smallskip

To discuss the proposed questions in further detail, we build on a convenient description of IWR lattices which we outline next. An important equivalence relation on lattices is geometric similarity: two lattices $\Lambda_1, \Lambda_2 \subset \real^N$ are called {\it similar}, denoted $\Lambda_1 \sim \Lambda_2$, if there exists a positive real number $\alpha$ and an $N \times N$ real orthogonal matrix $U$ such that $\Lambda_2 = \alpha U \Lambda_1$. It is easy to see that similar lattices have the same algebraic  structure, i.e., for every sublattice $\Gamma_1$ of a fixed index in $\Lambda_1$ there is a sublattice $\Gamma_2$ of the same index in $\Lambda_2$ so that $\Gamma_1 \sim \Gamma_2$. Most geometric and optimization properties of lattices (such as packing density, covering thickness, kissing number, signal-to-noise ratio, etc.) are invariant on similarity classes. Moreover, a WR lattice can only be similar to another WR lattice, so it makes sense to speak of WR similarity classes of lattices. If $\Lambda \subset \real^2$ is a full rank WR lattice, then its set of minimal vectors $S(\Lambda)$ contains 4 or 6 vectors, and this number is 6 if and only if $\Lambda$ is similar to the hexagonal lattice 
$$\H := \begin{pmatrix} 2 & 1 \\ 0 & \sqrt{3} \end{pmatrix} \zed^2$$
(see, for instance Lemma 2.1 of \cite{lf:petersen}). Any two linearly independent vectors $\bx,\bwy \in S(\Lambda)$ form a minimal basis. While this choice is not unique, it is always possible to select $\bx,\bwy$ so that the angle $\theta$ between these two vectors lies in the interval $[\pi/3,\pi/2]$, and any value of the angle in this interval is possible. From now on when we talk about a minimal basis for a WR lattice in the plane, we will always mean such a choice. Then the angle between minimal basis vectors is an invariant of the lattice, and we call it the {\it angle of the lattice} $\Lambda$, denoted $\theta(\Lambda)$; in other words, if $\bx,\bwy$ is any minimal basis for $\Lambda$ and $\theta$ is the angle between $\bx$ and $\bwy$, then $\theta = \theta(\Lambda)$ (see \cite{hex} for details and proofs of the basic properties of WR lattices in $\real^2$). In fact, it is easy to notice that two WR lattices $\Lambda_1,\Lambda_2 \subset \real^2$ are similar if and only if $\theta(\Lambda_1)=\theta(\Lambda_2)$ (see \cite{hex} for a proof). Therefore the set of all similarity classes of WR lattices in $\real^2$ is bijectively parameterized by the set of all possible values of the angle, which is the interval $[\pi/3,\pi/2]$. On the other hand, this parameterization becomes less trivial if we talk about similarity classes of planar IWR lattices. In other words, one may wonder what are the possible values of $\theta(\Lambda)$ in the interval $[\pi/3,\pi/2]$ if $\Lambda$ is IWR?

The following parameterization follows from the classical theory of integral lattices and quadratic forms (see, for instance Chapter 1 of~\cite{martinet}). 

\begin{prop} \label{param} Let $\Lambda \subset \real^2$ be an IWR lattice, then
\begin{equation}
\label{cos_sin}
\cos \theta(\Lambda) = \frac{p}{q},\ \sin \theta(\Lambda) = \frac{r\sqrt{D}}{q}
\end{equation}
for some $p,r,q,D \in \zed_{>0}$ such that
\begin{equation}
\label{prqD}
p^2+Dr^2=q^2,\ \gcd(p,q)=1,\ \frac{p}{q} \leq \frac{1}{2}, \text{ and } D \text{ squarefree},
\end{equation}
and so $\Lambda$ is similar to
\begin{equation}
\label{OprqD}
\Omega_D(p,q) := \begin{pmatrix} q & p \\ 0 & r\sqrt{D} \end{pmatrix} \zed^2.
\end{equation}
Moreover, for every $p,r,q,D$ satisfying \eqref{prqD}, $\Omega_D(p,q)$ is an IWR lattice with the angle $\theta(\Omega_D(p,q))$ satisfying \eqref{cos_sin}, and $\Omega_D(p,q) \sim \Omega_{D'}(p',q')$ if and only if $(p,r,q,D) = (p',r',q',D')$. In addition, if $\Lambda$ is any IWR lattice similar to $\Omega_D(p,q)$, then
\begin{equation}
\label{min_lattice}
\left| \Lambda \right| \geq \left| \frac{1}{\sqrt{q}} \Omega_D(p,q) \right|,
\end{equation}
where the lattice $\frac{1}{\sqrt{q}} \Omega_D(p,q)$ is also IWR. Due to this property, we call $\frac{1}{\sqrt{q}} \Omega_D(p,q)$ a {\it minimal} IWR lattice in its similarity class.
\end{prop}

\begin{rem} Notice in particular that the integer lattice $\zed^2 = \Omega_1(1,1)$ and the hexagonal lattice $\H = \Omega_3(1,2)$.
\end{rem}

Hence we see that the set of similarity classes of planar IWR lattices is in bijective correspondence with the set of 4-tuples $(p,r,q,D)$ satisfying \eqref{prqD}. Here is an explicit characterization of this set of 4-tuples, which will be useful to us.

\begin{lem} \label{pqrD_par} Let $D$ be a positive squarefree integer and $m,n \in \zed$ with $\gcd(m,n)=1$ and $\sqrt{\frac{D}{3}} \leq \frac{m}{n} \leq \sqrt{3D}$. Now $p,r,q,D \in \zed_{>0}$ satisfy \eqref{prqD} if and only if
\begin{equation}
\label{mn_par}
p =  \frac{| m^2-Dn^2 |}{2^e \gcd(m,D)},\ r = \frac{2mn}{2^e \gcd(m,D)},\ q = \frac{m^2 + Dn^2}{2^e \gcd(m,D)},
\end{equation}
where
\begin{equation}
\label{e_def}
e = \left\{ \begin{array}{ll}
0 & \mbox{if either $2 \mid D$, or $2 \mid (D+1), mn$} \\
1 & \mbox{otherwise.}
\end{array}
\right.
\end{equation}
\end{lem}

We prove Lemma \ref{pqrD_par} in Section~\ref{parameter}. Further, let us say that an IWR planar lattice $\Lambda$ is of {\it type $D$} for a squarefree $D \in \zed_{>0}$ if it is similar to some $\Omega_D(p,q)$ as in \eqref{OprqD}. The type is uniquely defined, i.e., $\Lambda$ cannot be of two different types. Moreover, a planar IWR lattice $\Lambda$ is of type $D$ for some squarefree $D \in \zed_{>0}$ if and only if all of its IWR finite index sublattices are also of type $D$. If this is the case, $\Lambda$ contains a sublattice similar to $\Omega_D(p,q)$ for every 4-tuple $(p,r,q,D)$ as in~\eqref{prqD}. Hence the set of planar IWR lattices is split into types which are indexed by positive squarefree integers with similarity classes inside of each type $D$ being in bijective correspondence with solutions to the ternary Diophantine equation $p^2+r^2D=q^2$ as parameterized in Lemma~\ref{pqrD_par}. 

\begin{rem} \label{composition} In fact, the set of similarity classes of IWR lattices of a fixed type can be endowed with a semigroup structure, coming from the geometric group law on rational points of a Pell conic; we include a brief discussion of this fact in Section~\ref{parameter} below (Lemma~\ref{conic}). The correspondence between IWR lattices and solutions to the Pell-type equations as described above follows from the theory of integral quadratic forms, as we indicated; it can also be obtained by an elementary argument, however we do not include it here in the interest of brevity of the exposition.
\end{rem}
\smallskip

In Section~\ref{optim} we discuss a possible connection between Questions~\ref{quest2} and~\ref{quest3}, and then use the above-described correspondence to provide an algorithmic procedure in answer to Question~\ref{quest2}.

\begin{thm} \label{optimize} A positive real number $\Delta$ is a determinant value of IWR lattices if and only if $\Delta = M\sqrt{D}$ where $M,D \in \zed_{>0}$ with $D$ squarefree so that the set
\begin{equation}
\label{mn_set}
\mn(M) = \left\{ (m,n) \in \zed^2_{>0} : \gcd(m,n)=1, \sqrt{\frac{D}{3}} \leq \frac{m}{n} \leq \sqrt{3D}, \frac{2mn}{2^e \gcd(m,D)}\ \Big|\ M \right\}
\end{equation}
where $e$ is as in \eqref{e_def}, is not empty. Fix such a $\Delta$, and let $(m,n) \in \mn(M)$ be the pair that maximizes the expression
$$\frac{m}{n} + D \frac{n}{m}$$
on $\mn(M)$. Now define $p,r,q$ as in \eqref{mn_par} for this choice of $m,n$ and let $k=M/r$. Then 
\begin{equation}
\label{max_SNR}
\Lambda = \sqrt{\frac{k}{q}}\ \Omega_D(p,q)
\end{equation}
is an IWR lattice with $\det(\Lambda)=\Delta$ and $|\Lambda|=kq$ which maximizes $|\Lambda|$  among all planar IWR lattices with determinant $\Delta$. This lattice can be found in a finite number of steps for each fixed $\Delta$.
\end{thm}

\begin{rem} \label{ex_table} Some examples of such norm-maximizing lattices are presented in Table~\ref{table1} below.
\end{rem}

In Section~\ref{cnt} we obtain the following counting estimate, which answers Question~\ref{quest1}.

\begin{thm} \label{count} For $\Delta \in \real_{>0}$, define $\IWR(\Delta)$ to be the set of all planar IWR lattices, up to rotation and reflection, with determinant $=\Delta$. Then the set $\IWR(\Delta)$ is finite for any $\Delta$, and it is only nonempty if $\Delta=M\sqrt{D}$ with the set $\mn(M)$ as in \eqref{mn_set} nonempty. In this latter case, the cardinality of the set $\IWR(\Delta)$ satisfies
\begin{equation}
\label{IWR_est}
\left| \IWR(\Delta) \right| \leq \frac{1}{2} \sum_{r \mid M} 2^{\omega(rD)}.
\end{equation}
Moreover, 
\begin{equation}
\label{IWR_est0}
\left| \IWR(\Delta) \right| \ll \sum_{r \mid M} \sum_{g \mid r} \mu \left( \frac{r}{g} \right) \frac{\tau(g^2D)}{\sqrt{\omega(gD)}},
\end{equation}
where $\tau(u)$ is the number of divisors, $\omega(u)$ is the number of prime divisors, and $\mu(u)$ is the M\"obius function of an integer $u$. The constant in the Vinogradov notation $\ll$ does not depend on $\Delta$.
\end{thm}

We are now ready to proceed.
\bigskip

\section{Parameterization lemmas}
\label{parameter}

In this section we start by proving Lemma~\ref{pqrD_par}. The following lemma is used in the proof, which we state here for the reader's convenience.

\begin{lem} [Lemma 2.1 of \cite{wr3}] \label{dioph} Consider the Diophantine equation
\begin{equation}
\label{d1}
\alpha x^2 + \beta xy + \gamma y^2 = \delta z^2,
\end{equation}
where $\alpha,\beta,\gamma,\delta \in \zed$ with $\beta^2 \neq 4\alpha\gamma$ and $\delta \neq 0$. Then either this equation has no integral solutions with $z \neq 0$, or all such solutions $(x,y,z)$ of \eqref{d1} are rational multiples of
\begin{equation}
\label{sol}
\begin{split}
x & = \gamma n (an - 2bm) - (\alpha a + \beta b) m^2, \\
y & = \alpha m (bm - 2an) - (\gamma b + \beta a) n^2, \\ 
z & = \pm c (\alpha m^2 + \beta mn + \gamma n^2),
\end{split}
\end{equation}
where $m,n \in \zed$ with $\gcd(m,n)=1$ and $m \geq 0$; here $(a,b,c)$ is any integral solution to \eqref{d1} with $c \neq 0$. In this later case, every multiple of \eqref{sol} is a solution to \eqref{d1} by homogeneity of the equation \eqref{d1}.
\end{lem}

\proof[Proof of Lemma \ref{pqrD_par}]
We start by applying Lemma~\ref{dioph} to the equation $p^2+Dr^2=q^2$ for a fixed squarefree $D$: since $(p,r,q)=(1,0,1)$ is an integral solution of this equation with $q \neq 0$, the lemma guarantees that all positive integral solutions of this equation with $q \neq 0$ are rational multiples of
\begin{equation}
\label{mn_1}
p_0=|m^2-Dn^2|,\ r_0=2mn,\ q_0=m^2+Dn^2,
\end{equation}
where $m,n$ range over all relatively prime non-negative integers, not both 0. In order for $(p,r,q,D)$ to satisfy \eqref{prqD}, we need two more conditions: $\gcd(p,q)=1$ and $p/q \leq 1/2$. First consider $p_0,r_0,q_0$ as in \eqref{mn_1} and notice that the fact that $p_0^2+Dr_0^2=q_0^2$ implies that $\gcd(p_0,q_0)=\gcd(r_0,q_0)=\gcd(p_0,r_0,q_0)$. Since $\gcd(m,n)=1$, it is easy to notice that $\gcd(p_0,q_0) = 2^e \gcd(m,D)$, where $e$ is as in \eqref{e_def}. Hence if we define $p,r,q$ as in \eqref{mn_par}, we ensure that they are relatively prime, and this covers all the relatively prime solutions of our equation for each fixed $D$. Finally, we need to select only the solutions with $p/q \leq 1/2$, which means that
$$-1/2 \leq \frac{m^2-Dn^2}{m^2+Dn^2} \leq 1/2,$$
and so we must have
\begin{equation}
\label{mn_ineq}
\sqrt{\frac{D}{3}} \leq \frac{m}{n} \leq \sqrt{3D}.
\end{equation}
This completes the proof of the theorem.
\endproof

We also briefly mention the algebraic structure of the planar IWR lattices.

\begin{lem} \label{conic} Let $D>0$ be squarefree and let $\C(D)$ be the set of similarity classes of all IWR lattices of type $D$. Let us write $C_D(p,q)$ for each such class, i.e., for each $(p,q)$ satisfying \eqref{prqD},
\begin{equation}
\label{CD}
C_D(p,q) = \left\{ \Lambda : \Lambda \sim \Omega_D(p,q) \right\},
\end{equation}
and so
$$\C(D) = \{ C_D(p,q) : (p,q) \text{ satisfy \eqref{prqD}} \}.$$
Then the set $\C(D)$ has the structure of an abelian semigroup, induced by the composition law on rational points of the Pell conic corresponding to $D$.
\end{lem}

\proof
A Pell conic is a curve given by the equation $x^2-Dy^2=1$. The following commutative composition law on the set of rational points on a Pell conic is defined in \cite{conics}:
\begin{equation}
\label{con_+}
(x_1,y_1) + (x_2,y_2) = (x_1x_2+Dy_1y_2, x_1y_2+x_2y_1).
\end{equation}
In \cite{conics}, this operation is also described geometrically by analogy with addition on an elliptic curve. Notice that a rational point $(x,y)=(q/p,r/p)$ is on this curve if and only if
\begin{equation}
\label{con_prq}
p^2+r^2D=q^2.
\end{equation}
Then \eqref{con_+} induces the following commutative composition law on the set of solutions $(p,r,q)$ of \eqref{con_prq}:
\begin{equation}
\label{prq_+}
(p_1,r_1,q_1) + (p_2,r_2,q_2) = \frac{1}{g} (p_1p_2, r_1q_2+r_2q_1, q_1q_2+Dr_1r_2),
\end{equation}
where $g = \gcd(p_1p_2, r_1q_2+r_2q_1, q_1q_2+Dr_1r_2)$. It is easy to check that the set of solutions of \eqref{con_prq} is closed under this operation. Moreover, since $D>0$,
$$\frac{q_1q_2+Dr_1r_2}{p_1p_2} \geq \frac{q_1}{p_1} \times \frac{q_2}{p_2},$$
and so whenever $p_1/q_1, p_2/q_2 \leq 1/2$, we will have
$$\frac{p_1p_2}{q_1q_2+Dr_1r_2} \leq \frac{1}{4}.$$
This ensures that $\C(D)$ is closed under this operation, and hence has a structure of an abelian semigroup, although not a monoid: the point $(1,0,1)$, which serves as identity, is not in $\C(D)$.
\endproof
\bigskip

\section{Optimization properties}
\label{optim}

In this section we investigate Questions~\ref{quest2} and~\ref{quest3}. Let $\Lambda$ be a planar IWR lattice, then
\begin{equation}
\label{det_1}
\Lambda = \sqrt{\frac{k}{q}} U \Omega_D(p,q)
\end{equation}
for some $(p,r,q,D)$ as in \eqref{prqD}, $k \in \zed_{>0}$, and a $2 \times 2$ real orthogonal matrix $U$. Now suppose that $\Delta = M\sqrt{D}$, $M \in \zed_{>0}$, is fixed and let $\Lambda \in \IWR(\Delta)$ be given as in \eqref{det_1} so that $kr=M$. Then 
\begin{equation}
\label{min_norm}
|\Lambda| = kq = \frac{Mq}{r},
\end{equation}
and so to maximize $|\Lambda|$ on $\IWR(\Delta)$ we need to maximize $q/r$. A trivial upper bound for $|\Lambda|$ is given by $\frac{2\Delta}{\sqrt{3}}$: this is just a restatement of the fact that $\Delta = |\Lambda| \sin \theta(\Lambda)$ and $\theta \in [\pi/3,\pi/2]$.
\smallskip

We start by discussing a connection between the problems of maximizing $|\Lambda|$ and minimizing $E_{\Lambda}(s)$ on sets of WR lattices of fixed determinant in~$\real^2$. This discussion is an  adaptation and correction of Lemma~5.2 of \cite{wr3}.

\begin{lem} \label{min_epstein} Let $\Delta$ be a positive real number, and let $\WR_2(\Delta)$ be the set of all full rank WR lattices in $\real^2$ with determinant $\Delta$. Then for any fixed real number $s \geq 3$, $E_{\Lambda}(s)$ is a decreasing function of $|\Lambda|$ on $\WR_2(\Delta)$.
\end{lem}

\proof
Let $Q_{\Lambda}(x,y)$ be the quadratic form of $\Lambda$ corresponding to a minimal basis, then
\begin{equation}
\label{QL}
Q_{\Lambda}(x,y) = |\Lambda| (x^2+y^2+2xy\cos \theta),
\end{equation}
where $\theta = \theta(\Lambda) \in [\pi/3,\pi/2]$ and $|\Lambda|$ is as in \eqref{min_norm}. Now 
\begin{equation}
\label{cos_T}
\cos \theta = \frac{\sqrt{|\Lambda|^2 - \Delta^2}}{|\Lambda|} = \sqrt{ 1 - \frac{\Delta^2}{|\Lambda|^2}},
\end{equation}
and $0 \leq \cos \theta \leq 1/2$. Lemma~1 of \cite{ryskov} guarantees that $E_{\Lambda}(s)$ is a decreasing function of $\cos \theta$ for any real $s \geq 3$, and \eqref{cos_T} implies that $\cos \theta$ is an increasing function of $|\Lambda|$. Hence $E_{\Lambda}(s)$ is a decreasing function of $|\Lambda|$ on $\WR_2(\Delta)$ for $s \geq 3$.
\endproof

In fact, it seems likely that the statement of Lemma~\ref{min_epstein} should hold for smaller real values of $s$ as well. At the very least, we have the following bounds.

\begin{lem} \label{SNR_min} With notation as in Lemma~\ref{min_epstein}, let $s>1$ be real. Then there exist real constants $C_1(s)$ and $C_2(s)$, dependent only on $s$, such that
\begin{equation}
\label{E_bounds}
\frac{C_1(s)}{|\Lambda|^s} \leq E_{\Lambda}(s) \leq \frac{C_2(s)}{|\Lambda|^s},
\end{equation}
for every $\Lambda \in \WR_2(\Delta)$. 
\end{lem}

\proof
Combining \eqref{QL} and \eqref{cos_T}, we obtain
$$Q_{\Lambda}(x,y) = T x^2+ T y^2+2xy \sqrt{T^2-\Delta^2},$$
where $T=|\Lambda|$. The Epstein zeta-function of $\Lambda$ is then given by
\[
\begin{split}
E_{\Lambda}(s) = \sum_{x,y \in \zed \setminus \{0\}} Q_{\Lambda}(x,y)^{-s} = \sum_{x,y \in \zed \setminus \{0\}} \frac{1}{\left( T x^2+ T y^2+2xy \sqrt{T^2-\Delta^2} \right)^s} \\
 = \sum_{x,y \in \zed_{>0}} \left( \frac{2}{\left( T x^2+ T y^2+2xy \sqrt{T^2-\Delta^2} \right)^s} + \frac{2}{\left( T x^2+ T y^2-2xy \sqrt{T^2-\Delta^2} \right)^s} \right).
\end{split}
\]
Now recall that since $\theta \in [\pi/3,\pi/2]$, we must have $\frac{\sqrt{3}T}{2} \leq \Delta \leq T$, and so $0 \leq \sqrt{T^2-\Delta^2} \leq T/2$. Hence for each fixed real $s > 1$, we have
\begin{equation}
\label{E_bnd_1}
E_{\Lambda}(s) \leq \frac{2}{T^s} \sum_{x,y \in \zed_{>0}} \left( \frac{1}{\left( x^2+ y^2 \right)^s} + \frac{1}{\left( x^2+ y^2 - xy \right)^s} \right),
\end{equation}
and
\begin{equation}
\label{E_bnd_2}
E_{\Lambda}(s) \geq \frac{2}{T^s} \sum_{x,y \in \zed_{>0}} \left( \frac{1}{\left( x^2+ y^2 \right)^s} + \frac{1}{\left( x^2+ y^2 + xy \right)^s} \right).
\end{equation}
Since both series in the bounds of \eqref{E_bnd_1} and \eqref{E_bnd_2} converge, we have \eqref{E_bounds}.
\endproof

\begin{rem} \label{E_min_heuristics} Since WR lattice $\Lambda$ with fixed $|\Lambda|$ and $\det(\Lambda)$ is unique up to multiplication by an orthogonal matrix $U$ (which does not change the value of $E_{\Lambda}(s)$ for any $s$), Lemmas~\ref{min_epstein} and~\ref{SNR_min} make it natural to expect that the total interference of $\Lambda$ is minimized on $\WR_2(\Delta)$ (and so $\SNR(\Lambda)$ is maximized) if and only if $|\Lambda|$ is maximized. 
\end{rem}

We are now ready to answer Question~\ref{quest2}.

\proof[Proof of Theorem~\ref{optimize}]
We will now discuss a finite procedure to maximize $q/r$, and hence $|\Lambda|$, on the set $\IWR(\Delta)$ using finiteness of this set along with Lemma~\ref{pqrD_par}. First we notice that $r$ has to be a divisor of $M=\Delta/\sqrt{D}$, hence we can start by going through the list of all possible divisors of $M$. For each such divisor $r$, consider all possible decompositions
$$r = \frac{2mn}{2^e \gcd(m,D)}$$
with relatively prime $m,n$ so that $m/n \in \left[ \sqrt{D/3}, \sqrt{3D} \right]$, as in \eqref{mn_par}.  Out of all such decompositions, we want to pick one which maximizes the ratio
$$q/r = \frac{m^2+Dn^2}{2mn} = \frac{1}{2} \left( \frac{m}{n} + D \frac{n}{m} \right).$$
This can be done in a finite number of steps, since there are finitely many values for $r$, a divisor of $M$, and for each $r$ there are finitely many such $m,n$. Hence we can choose $\Lambda$ maximizing $|\Lambda|$ and $\SNR(\Lambda)$ on $\IWR(\Delta)$ to be as in \eqref{max_SNR}. In particular, our argument confirms that $\Delta$ is a determinant value of an IWR lattice if and only if it is of the form $M\sqrt{D}$ with the set $\mn(M)$ as in \eqref{mn_set} nonempty. This completes the proof.
\endproof

\begin{rem} \label{mnx}
Let us write $m/n = \sqrt{D}x$ for appropriate $x \in \left[ 1/\sqrt{3}, \sqrt{3} \right]$. Then
$$q/r = \frac{\sqrt{D}}{2} (x + 1/x).$$
Now the function $f(x) = x+1/x$ assumes its maximal values on the interval $\left[ 1/\sqrt{3}, \sqrt{3} \right] $ at the endpoints and has a minimum at $x=1$. Hence, to maximize $q/r$ one should consider $m,n$ with $m/n$ close to the endpoints of the interval $\left[ \sqrt{D/3}, \sqrt{3D} \right]$. Keeping these considerations in mind can reduce the number of computational steps necessary to find maximizer for $|\Lambda|$ in $\IWR(\Delta)$ in each particular case.
\end{rem}

We give some computational examples in Table~\ref{table1} below.

\begin{center} 
\begin{table}[!ht]
\caption{Examples of IWR lattices $\Lambda$ with  $\det(\Lambda)=\Delta$ that maximize $|\Lambda|$ on $\IWR(\Delta)$} 
\begin{tabular}{|l|l|l|l|} \hline
$\Delta$&$|\Lambda|$&$\Lambda$ \\
\hline
$24\sqrt{5}$ & $61$ & $\sqrt{\frac{1}{61}}\Omega_5(29,61)$\\ \hline
$24\sqrt{7}$ & $69$ & $\sqrt{\frac{3}{23}}\Omega_7(9,23)$\\ \hline
$20\sqrt{11}$ & $75$ & $\sqrt{\frac{1}{3}}\Omega_{11}(7,15)$\\ \hline
$24\sqrt{13}$ & $98$ & $\sqrt{\frac{2}{49}}\Omega_{13}(7,15)$\\ \hline
$24\sqrt{17}$ & $104$ & $\sqrt{\frac{8}{13}}\Omega_{17}(4,13)$\\ \hline
$105\sqrt{19}$ & $510$ & $\sqrt{\frac{15}{34}}\Omega_{19}(15,34)$\\ \hline
$96\sqrt{23}$ & $522$ & $\sqrt{\frac{6}{87}}\Omega_{23}(41,87)$ \\ \hline
\end{tabular}
\label{table1}
\end{table}
\end{center}
\bigskip

\section{Counting estimates}
\label{cnt}

Here we prove the counting estimates of Theorem~\ref{count}.

\proof[Proof of Theorem~\ref{count}]
Let $\Delta=M\sqrt{D}$, as above, so that $\mn(M)$ defined in \eqref{mn_set} is nonempty. First we observe that such sets are in fact finite up to rotation and reflection - this is an immediate consequence of a more general fact that there are only finitely many isometry classes of integral lattices of fixed determinant in a fixed dimension (see remarks on p.~432 of~\cite{martinet}). Suppose now that $\Lambda \in \IWR(\Delta)$, then we can assume without loss of generality that
$$\Lambda = \sqrt{\frac{k}{q}}\ \Omega_D(p,q),$$
where $k=M/r$ and $p,r,q$ are as in \eqref{mn_par} for some $(m,n) \in \mn(M)$. Hence the choice of $p,r,q$ determines $\Lambda$ uniquely. For each $r \mid M$ define
\begin{equation}
\label{r_count}
f(r) = \left| \left\{ (p,q) \in \zed_{>0}^2 : q^2-p^2 = r^2D,\ \gcd(p,q) = 1,\ 0< \frac{p}{q} \leq \frac{1}{2} \right\} \right|,
\end{equation}
then
$$\left| \IWR(\Delta) \right| = \sum_{r \mid M} f(r).$$
Hence we want to produce estimates on $f(r)$. Define
$$f_1(r) = \left| \left\{ (p,q) \in \zed_{>0}^2 : q^2-p^2 = r^2D,\ \gcd(p,q) = 1 \right\} \right|,$$
and
\begin{equation}
\label{f_21}
f_2(r) = \left| \left\{ (p,q) \in \zed_{>0}^2 : q^2-p^2 = r^2D,\ 0< \frac{p}{q} \leq \frac{1}{2} \right\} \right|,
\end{equation}
and notice that
\begin{equation}
\label{f_f1_f2}
f(r) \leq \min \{ f_1(r), f_2(r) \},
\end{equation}
meaning that
\begin{equation}
\label{IWR_count}
\left| \IWR(\Delta) \right| \leq \sum_{r \mid M} \min \{ f_1(r), f_2(r) \}.
\end{equation}
The function $f_1(r)$ is well-studied; in particular, the following formula follows from Theorem 6.2.4 of \cite{mollin}:
\begin{equation}
\label{f1_est}
f_1(r) = \left\{ \begin{array}{ll}
2^{\omega(r^2D)-1} & \mbox{if $2 \nmid r^2D$, $r^2D > 1$} \\
2^{\omega(r^2D)-1} & \mbox{if $8 \mid r^2D$, $r^2D$ has odd prime divisors} \\
1 & \mbox{if $r^2D$ is a power of 2} \\
0 & \mbox{otherwise,}
\end{array}
\right.
\end{equation}
hence $f_1(r) \leq 2^{\omega(r^2D)-1} = 2^{\omega(rD)-1}$. Now \eqref{IWR_est} follows upon combining \eqref{IWR_count}, \eqref{f1_est}.

Next we estimate $f_2(r)$. Let $c=r^2D$, and let us write
\begin{equation}
\label{ab}
a=q-p,\ b=q+p,
\end{equation}
then $q=(a+b)/2$, $p=(b-a)/2$, and $ab=c$. Let $\alpha := p/q$, and assume that $0 < \alpha \leq 1/2$. Then let $\nu =\frac{1+\alpha}{1-\alpha}$, and observe that
$$1 <  \nu = \frac{b}{a} \leq 3.$$
Since $ab=c$, we have $b = \sqrt{\nu c}$, and so
$$\sqrt{c} < b \leq \sqrt{3c}.$$
Therefore
\begin{equation}
\label{f_22}
f_2(r) = \left| \left\{ b \in \zed_{>0} : b \mid c,\ \sqrt{c} < b \leq \sqrt{3c} \right\} \right|.
\end{equation}
For a positive integer $t$, Hooley's $\HDelta$-function of $t$ (see \cite{hall} for detailed information) is defined as
$$\HDelta(t) = \max_x \left| \left\{ b \in \zed_{>0} : b \mid t,\ e^x < b \leq e^{x+1} \right\} \right|.$$
Take $x=\log \sqrt{c}$, then
$$\left\{ b \in \zed_{>0} : b \mid c,\ \sqrt{c} < b \leq \sqrt{3c} \right\} \subseteq \left\{ b \in \zed_{>0} : b \mid c,\ e^x < b \leq e^{x+1} \right\},$$
since $\sqrt{3} < e$, and so $f_2(r) \leq \HDelta(c)$. Therefore an estimate on $f_2(r)$ would follow from estimates on $\HDelta(c)$, some of which can be found in Section~2 of \cite{wr1}; in particular, equations (10)-(13) of \cite{wr1} imply that the bound
\begin{equation}
\label{f2_est}
f_2(r) \leq O \left( \frac{\tau(c)}{\sqrt{\omega(c)}} \right) \leq O \left( c^{\frac{(1+\eps) \log 2}{\log \log c}} \right)
\end{equation}
holds for any $\eps > 0$, assuming $c$ is greater than some $c_0(\eps)$ for the second inequality; here  the constant in $O$-notation is independent of $c$.

Next notice that if $q^2-p^2=r^2D$ and $g \mid p,q$, then $g \mid r$, since $D$ is squarefree. This implies that
\begin{equation}
\label{f2r_mob}
f_2(r) = \sum_{g \mid r} f \left(\frac{r}{g} \right).
\end{equation}
Recall that the M\"obius function is defined by
$$\mu(u) = \left\{ \begin{array}{ll}
(-1)^{\omega(u)} & \mbox{if $u$ is squarefree} \\
0 & \mbox{otherwise,}
\end{array}
\right.$$
then applying the M\"obius inversion formula to \eqref{f2r_mob}, we obtain
\begin{equation}
\label{f_r_g}
f(r) = \sum_{g \mid r} \mu \left( \frac{r}{g} \right) f_2(g) \ll \sum_{g \mid r} \mu \left( \frac{r}{g} \right) \frac{\tau(g^2D)}{\sqrt{\omega(g^2D)}},
\end{equation}
by \eqref{f2_est}. This establishes \eqref{IWR_est0} upon the observation that $\omega(g^2D)=\omega(gD)$.
\endproof

\begin{rem} \label{count_bound} Theorems 431 and 432 of \cite{hardy} state that normal orders of $\omega(u)$ and $\tau(u)$ are $\log \log u$ and $2^{\log \log u}$, respectively. This implies that one would normally expect
$$\frac{\tau(u)}{\sqrt{\omega(u)}} \leq 2^{\omega(u)}$$
for a randomly chosen integer $u$ (in the appropriate sense).
\end{rem}

\bigskip

{\bf Acknowledgment.} We would like to thank the Fletcher Jones Foundation-supported Claremont Colleges research experience program, under the auspices of which this work was done during the Summer of 2011. We thank Professor Wai Kiu Chan for his useful comments on the subject of this paper. We are also grateful to the anonymous referees for their many helpful suggestions.
\bigskip

\bibliographystyle{plain}  
\bibliography{iwr_lattices_DCG}    
\end{document}